\documentclass[11pt]{amsart}
\usepackage{amsmath,amssymb,enumerate}
\usepackage{epsf,epsfig,amsfonts,graphicx,color}  
 \usepackage[applemac]{inputenc} % (pour les accents)

 \date{}   %%%%%%% Set the correct date
  
 \numberwithin{equation}{section}   
\newtheorem{Th}{\rm\bf Theorem}

\newtheorem{remark}{\rm\bf Remark}
\theoremstyle{definition}

\newcommand{\weg}[1]{}
\newcommand{\be}{\begin{equation}}
\newcommand{\ee}{\end{equation}}

\title[Can we  make a Finsler metric complete?]{Can we make a Finsler metric complete by a trivial projective change?   } \date{}
\author{Vladimir S. Matveev }\thanks{Institute of Mathematics,  Friedrich-Schiller-Universit\"at Jena,  07737 Jena Germany\\  vladimir.s.matveev@gmail.com}

\begin{document}
%\keywords{Finsler metrics, geodesic equivalence, projective reparameterization} 

%	\subjclass{58B20, 53C60, 53C22, 53B10, 53A20} 
\begin{abstract} 
A trivial projective change of a Finsler metric $F$ is the Finsler metric $F + df$. I explain when it is possible to make a given Finsler metric both forward and backward complete by a trivial projective change. 

The problem actually came from lorentz geometry and  mathematical  relativity: it was observed that it is possible to understand  the light-line geodesics of  a (normalized, standard) stationary   4-dimensional  space-time  as geodesics of a certain Finsler Randers metric on a 3-dimensional manifold.     The trivial projective change of the Finsler metric corresponds to the  choice of another  3-dimensional slice, and the  existence of a  trivial projective change that is forward and backward complete is equivalent to the  global  hyperbolicity of the space-time. 
\end{abstract} 
\maketitle

\section{Statement of the problem, motivation and the main result}
Let $(M, F) $ be a connected Finsler manifold and $f:M\to \mathbb{R}$ be a function such that 
\begin{equation} \label{-1} 
F(x,v) + d_xf(v)>0  \  \textrm{  for all $(x,v)\in TM$ with $v\ne 0$}.\end{equation} By a {\it trivial projective change}
 we   	 understand the Finsler metric $F + df$.

It is customary  in Finsler geometry to require the Finsler metric to be  \emph{strongly convex}, 
that is the Hessian of the restriction of $F^2$ to $T_xM\setminus \{0 \}$ is assumed  to be positive definite for any 
$x \in M$.  Our results do not require this assumption and are valid also 
for Finsler metrics that are not strictly convex. Let us note though that if the metric $F$ is strictly   convex then  the trivial projective change $F + df$ is also strictly convex.

 The metric $F+ df$ has the same unparameterized geodesics as  $F$. Indeed, a forward-geodesic connecting 
 two points $x,y\in M$ is an extremal of the forward-length functional 
\begin{equation} 
\label{1}
 L^+_F(c):=  \int_{a}^b F(c(t), \dot c(t)) dt\end{equation} 
  in the set of all smooth curves $c:[a,b]\to M$ connecting $x$ and $y$. Now, replacing  $F$ by $F +df$ in \eqref{1}, we obtain 
  
  \begin{equation} \label{2}
  L^+_{F+df} (c):=  \int_{a}^b \left[F(c(t), \dot c(t)) +  d_{c(t)}f(\dot c  (t))  \right]dt= L^+_F(c)+ f(y)- f(x). \end{equation}
We see that the difference  $L^+_{F+df}(c) -L^+_F(c)$  is the constant $f(y)- f(x)$ so the extremals  of the   functional  \eqref{1}   are extremals of the functional \eqref{2}  and vice versa. 
   
Analogically, for the backward-length 
$$
L_F^-:= \int_{a}^b F(c(t), - \dot c(t)) dt, 
$$ we obtain  $L^-_{F+df}(c) -L^-_F(c)=f(x)- f(y)$ implying the backward-geodesics of $F$ and $F+ df $ 
coincide. 

Note that, though the unparameterized geodesics of $F$ and $F+ df $ coincide, the arc-length parameter of the geodesics and also the distance functions generated by the Finsler metrics do not coincide (unless $f$ is constant.) More precisely,  the forward and backward 
 distance functions 
 $$
 dist_{F}^\pm(x,y)  = \inf\{L_F^\pm(c) \mid c:[a,b]\to M \textrm{  with $c(a)= x$ , $c(b)=y$}\}$$
 and the corresponding distance functions    for $F+ df $ are related by 
 
\begin{equation} \label{3}
 dist_{(F+ df )}^\pm(x,y)= dist_F^\pm(x,y) \pm (f(y)- f(x)). 
 \end{equation}

 The goal of this note is to answer the questions under what conditions one could make a Finsler metric simultaneously forward  complete and backward complete  by an appropriate  trivial projective change. We will assume that all objects we consider in our paper are sufficiently smooth. The assumption that the  metric $F$ is smooth is very natural in view of our motivation (see \S \ref{motivation} below).  We  will see that the restriction 
  that the (searched)  function $f$ is smooth (which is also natural in view  of the motivation) actually makes our proof more complicated: if we allow the Lipschitz functions (and in Remark \ref{rem:1} we explain why we may  do it), the proof becomes  shorter and does  not require the Appendix where we proved that it is possible, for any $\varepsilon_1, \varepsilon_2>0$,  to $\varepsilon_1$- approximate an $1$-Lipschitz function by a $1+ \varepsilon_2$-Lipschitz function, where the Lipschitz property is understood with respect to the (nonsymmetric) distance function coming from the metric $F$.

\subsection{Motivation}  \label{motivation} 
Our  motivation to study this question came from the  mathematical relativity and Lorentz differential 
  geometry.   Following \cite{1,2,2a,2b}, see also references therein, we consider the  (normalized, standard) stationary space-time   $(M^4= \mathbb{R}\times S^3, G)$. Here    
  $S$ is a 3-dimensional manifold. The condition that the space-time is normalized, standard, stationary means that, in any   local coordinate system    $(t, x^1,x^2, x^3)$ where $t$ is the coordinate on $\mathbb{R}$ and $x^1,x^2,x^3$ are local coordinates on a 3-manifold $S$,  the  metric 
  $G$  is given by the formula 
\begin{equation}\begin{array}{ccl} G & = & -dt^2 + 2 \omega_i dx^i dt  + g_{ij} dx^idx^j\\ &=& -(dt - \omega_i dx ^i)(dt- \omega_j dx^j) + (g_{ij} + \omega_i \omega_j)dx^i dx^j  ,\end{array} \label{lorenz} \end{equation}
  where $g= g(x)_{ij} $ $i,j= 1,...,3$,  is a Riemannian metric  on $S$ and $\omega= \omega(x)_i$ 
  is a 1-form on $S$. \weg{ such that $g_{ij} + \omega_i\omega_j$ is positive defined  at every point of $S$.}  
  
  \begin{remark} Note that the above definition of the  normalized, standard, stationary spacetime is not the usual    one. Usually,  one 
  defines  a standard stationary spacetime as the one which is causal (i.e. it does not admit closed causal curves) and which admits a complete time-like Killing vector field $K$.  By \cite{R1}, this is equivalent to the condition  that $M$ is isometric to a product $\mathbb{R}\times S$, where $S$ is some (appropiate) space-like hypersurface, and $K$, in the coordinate system corresponding to this decomposition,   is   the  vector field $\tfrac{\partial  }{ \partial t} $. Then,  the metric can be  written locally as in \eqref{lorenz}; it is easy to check that   $\omega_i$ and $g_{ij}$ could be viewed as  objects globally defined on $S$   since   for $\xi \in TS$ we have  
  $\omega(\xi) = G\left(\tfrac{\partial  }{ \partial t} ,\xi \right) = G\left(K,\xi \right)$ and $g$ is simply the  restriction of $G$ to  $S$. 
  
  \end{remark} 
  
  Next, on  $S$, we consider the (Randers) Finsler metric 
  \begin{equation} \label{randers} 
  F(x,\dot x) = \sqrt{(g(x)_{ij} + \omega(x)_i\omega(x)_j)\dot x^i \dot x^j}  + \omega(x)_{i}\dot x^i.
  \end{equation}   
  
  As it was observed and actively studied in  \cite{1,2, DPS,4}, this Finsler metric  and the initial Lorentz metric $G$ are closely related. In particular, for every  
   light-like geodesic  $\gamma(\tau)= (t(\tau),x^1(\tau), x^2(\tau), x^3(\tau))$ of $G$, its ``projection'' to $S$, i.e., 
   the curve $\tau\mapsto (x^1(\tau), x^2(\tau), x^3(\tau))$ on $S$  is a (probably, reparameterized)  geodesic  of the Finsler metric \eqref{randers}.  Moreover, the slice  \begin{equation} 
 \{0\} \times S = \{(0,x)\mid x\in S\}\label{slice} \end{equation}  
    is a Cauchy hypersurface      of $(M,  G)$  if and only if the metric $F$ is forward and backward complete, see \cite[Theorem 4.4]{1}. 
  
  Note that it is  possible to take another decomposition of $M$ in the product of $\mathbb{R} \times S'$ such that the metric $G$  written in the coordinates adapted to the new decomposition still has the form \eqref{lorenz}  with possibly different $g$ and $\omega$.  
    
   Indeed,  consider another local coordinate systems $(t', x^1 ,x^2, x^3)$ such that $t'= t + f(x^1, x^2, x^3)$ (and the  coordinates $x^1,x^2, x^3$  are the same). Physically, this choice of the coordinates 
   corresponds to the choice of another space-like slice:    by the  ``old''  slice  we understand   the 3-dimensional submanifold  \eqref{slice},   
    and by the new one we understand 
   $ \{(f(x), x)\mid x\in S\}$.

   In the new coordinates $   (t', x^1 ,x^2, x^3) $, in view of $dt'= dt - df$,   the metric $G$  reads

   \begin{equation} -\left(dt' - \underbrace{\left(\omega_i + \tfrac{\partial f}{\partial x^i}\right)}_{\omega'_i}  dx^i  \right)\left(dt' -\underbrace{\left(\omega_j + \tfrac{\partial f}{\partial x^j}\right)}_{\omega'_j} dx^j\right)+  \left(g_{ij} +\omega_i\omega_j\right) dx^idx^j \label{randers2}.\end{equation}
   
   We see  that the Finsler metric \eqref{randers} constructed by the metric  \eqref{randers2} is related 
   to the initial  metric \eqref{randers} constructed by  \eqref{lorenz}  by the formula 
   $F'= F+ df$, i.e., is the trivial projective change of the metric \eqref{randers}.   It is easy to check that  the  slice  $ \{(f(x), x)\mid x\in S\}$ is space-like if and only if $$\sqrt{(g_{ij} + \omega_i\omega_j)\dot x^i \dot x^j}  + \omega_{i}\dot x^i + \tfrac{\partial f}{\partial x^i} \dot x^i \ \ \textrm{is positive for all  $\dot x^i \ne 0$},$$ i.e., if and only if  $f$ satisfies the condition from the definition of the 
    the trivial projective change with respect to  the Finsler metric \eqref{randers}.

   Thus, the question we are study, i.e., the existence of a trivial projective change of a Finsler metric such that the result is  forward and backward complete is, in the special case when the metric  is the Randers metric 
   coming  from the Lorentz metric \eqref{lorenz} by  the formula \eqref{randers}, 
   equivalent to the existence of the function $f$ such that the corresponding slice $ \{(f(x), x)\mid x\in S\}$ is a Cauchy hypersurface. Note that 
    if such Cauchy hypersurface exists  then the space-time is  { gobally hyperbolic}, see \cite{R3,R5} for details,  
    and  
   global hyperbolicity is an important condition to be studied in any space-time.

  It appears though that  the special form of the metric $F$ suggested by the motivation does not make (our version of) the answer  simpler, so we give the answer for the general Finsler metrics. It seems   that even in the well studied  situation when the Finsler metric is a Randers one, i.e., in the situation suggested by relativity, our main result which is Theorem \ref{thm1} below is new; cf. \cite[Theorem 5.10]{1}. Actually, \cite[Theorem 5.10]{1} and our main Theorem \ref{thm1} restricted to the Randers metrics  are very similar: the difference is that in \cite[Theorem 5.10]{1} one (essentially) assumes that the function $D^+ + D^-$  (see below) 
  is proper for  all choices of the  point $p$  as the initial point and we require this for one point $p$ only.

\subsection{Main result} 
 We fix  an arbitrary point $p\in M$  and consider the functions $D^+, D^- :M\to \mathbb{R}$ given by 
  $$ D^\pm(x) := dist^\pm(p, x).$$ 
  
  \begin{Th} \label{thm1}  $F$ can be made forward and backward  complete by a trivial  projective change if and only if the function $D^+ + D^-$ is proper. 
  \end{Th} 

Recall that  a (continuous) function is {\it proper}, if the preimage of every  compact set is compact  or empty.
Since the function $D^+ + D^-$ is nonnegative and $D^+(p) + D^-(p)= 0$, the function   $D^+ + D^-$ is proper 
if and only if for every $c\in \mathbb{R}_{\ge 0}$ the set 
\begin{equation} \label{form}  
\{x \in M\mid D^+(x)+ D^-(x) \le c\}\end{equation} 
is compact.

\section{ Proof of  Theorem \ref{thm1}.} 
First observe that   if the function $D^+ + D^-$ is proper  then 
the function $\alpha_1 D^+ + \alpha_2 D^-$   is proper for arbitrary  positive numbers $\alpha_i$, and that if the   function $\alpha_1 D^+ + \alpha_2 D^-$   is proper for some   positive numbers $\alpha_i$ then the function $D^+ + D^-$ is proper. Indeed, the set 

\begin{equation} \label{form1} \{ x\in {M}\mid  \alpha_1 D^+(x)  + \alpha_2 D^-(x) \le c\}\end{equation}  
is a (evidently, closed) subset of 
 $$
 \{ x\in M \mid  D^+(x) + D^-(x) \le \frac{c}{\alpha_{\min}}  \}, \ \ 
 \textrm{where $\alpha_{\min}= \min(\alpha_1, \alpha_2)$}
 $$
 Then, if all the sets of the form  \eqref{form} are compact, all sets of the form \eqref{form1} are compact as well implying  the function $\alpha_1 D^+(x) + \alpha_2 D^-(x)$ is proper.

 Now, if the function $\alpha_1 D^+ + \alpha_2 D^-$ is proper,  the set 
\eqref{form} is  compact as a closed subset of    
$$\{ x\in M\mid  \alpha_1 D^+  + \alpha_2 D^-  \le \alpha_{\max} c \}, \ \ \textrm{ where $\alpha_{\max} = \max(\alpha_1,\alpha_2)$},$$  implying $D^+ + D^-$ is proper. 

We will now  show  Theorem \ref{thm1} in  the direction ``$\Longrightarrow$'': we show that if the function $D^+ + D^-$ is not proper,  
then no projective change $F + df$ is complete. Let $R>0$ be the number such that the set $$B_R:= \{x \in M \mid D^+(x) + D^-(x)\le R\}$$ is not compact.  For any  $x \in B_R$, we have $D^+(x)\le R$ and $D^-(x)\le R$. 
Then, in view of \eqref{3}, we have  
\begin{equation}\label{5} \begin{array}{ll} dist_{F +df}^+ (p,x) &= D^+(x) +   (f(x)- f(p))\le R+ (f(x)- f(p)) \\ 
dist_{F +df}^-(p,x) &= D^-(x)   + (f(p)- f(x)) \le R+ (f(p)- f(x)).\end{array}\end{equation} 
Since 
$dist_{F +df}^\pm (p,x)\ge 0$,  we obtain   that $-R\le f(p)-  f(x)\le R$. Then, the set 
$B_R $ lies in the set $\{x \in M\mid dist_{F+ df}^+(p, x)  + dist^-_{F+ df}(p,x)  \le  3 R \},  $ which, in its turn, lies in the set    
\begin{equation} 
\{x \in M\mid dist_{F+ df}^+(p, x)\le  3R \}   \cap \{dist^-_{F+ df}(p,x)  \le  3 R \}.   \label{4}\end{equation} Would the metric $F + df $ be forward and backward complete, the set \eqref{4} would be compact implying all its closed subsets are compact which contradicts our assumption that $B_R$ 
is not compact. Theorem is proved in one direction.

In order to prove it in the other direction,  we consider   the function $$f:M\to \mathbb{R}, \ \ f(x)   := \frac{D^-(x) - D^+(x)}{2}.$$

The function $f$ is not a priori smooth; 
next we   show that the function is ${1}{}$-Lipshitz w.r.t. the  distance 
$dist^+$, that is,  for every $x,y \in M$ we have \begin{equation} dist^+(x,y)\ge f(x)- f(y). \label{lip}\end{equation}
 
\begin{figure}
\includegraphics[width=.4\textwidth]{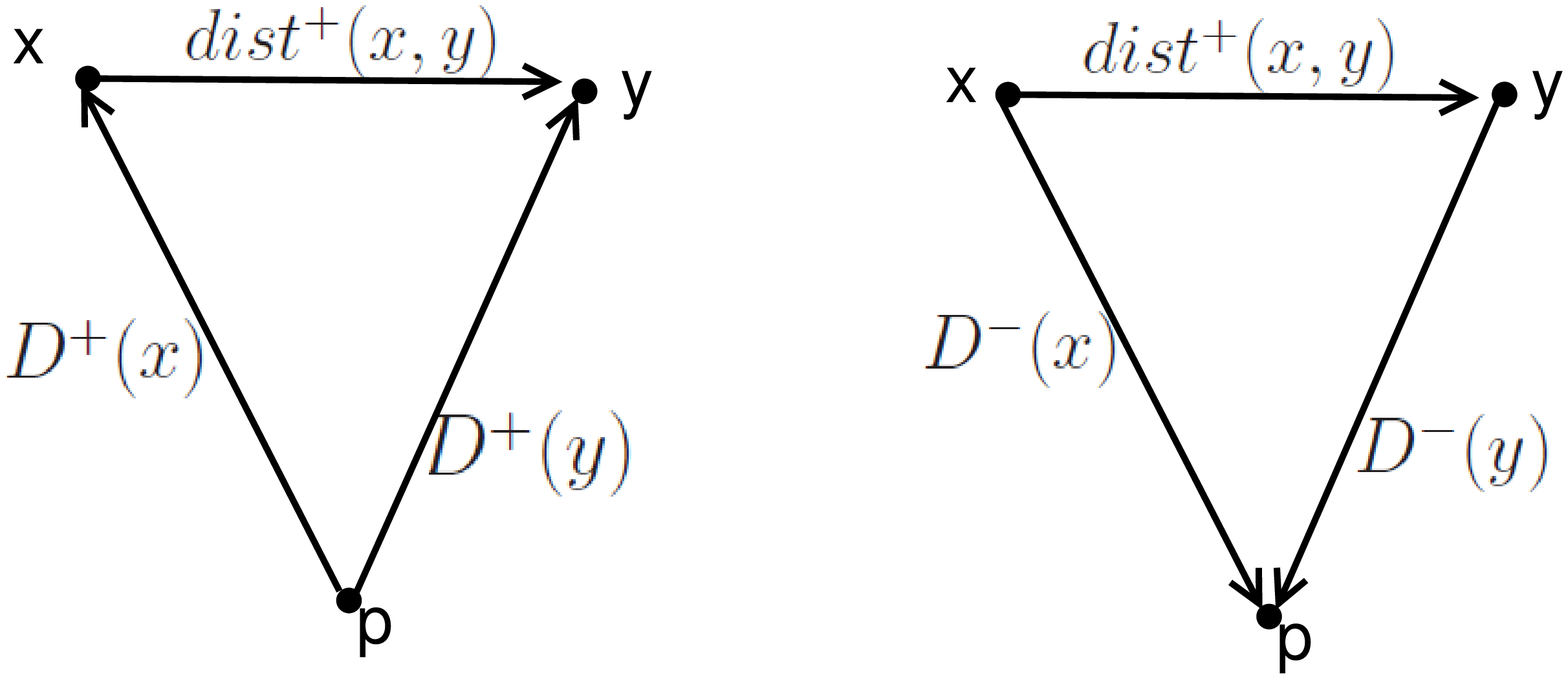}
\caption{Triangles for the triangle inequalities.} \label{fig.1}
\end{figure}

  Indeed, consider the triangles on the Fig. \ref{fig.1} and the corresponding triangle inequalities: 
\begin{equation} \label{10} D^+(x) +dist^+(x,y) \ge D^+ (y) \ , \ \  D^-(y)+ dist^+(x,y)\ge D^-(x). \end{equation}
The sum of these inequalities is equivalent to  \eqref{lip}

\begin{remark} 
As a consequence we obtain that the function $f$ is also 1-Lipschitz w.r.t. the symmetrized distance $dist^{sym}:= dist^+ + dist^-$. 
Since, locally, in a sufficiently small neighborhood, we can evidently find a euclidean structure such that the corresponding distance is not less than $dist^{sym}$, the function $f$ is  locally Lipschitz w.r.t. an Euclidean structure and is therefore differentiable at almost every point.  Moreover, the restriction of the function to every smooth curve is a locally lipschitz function and the formula \eqref{2} remains  valid though $df$ is not everywhere defined.  
\end{remark}

Let us now   take a  smooth function $\tilde f$ on $M$ such that 
\begin{enumerate} 
\item $|\tilde f(x) - f(x)| \le 1$  for all $x\in M$. 
\item $\tilde f(x)$ is 1.5-Lipschitz w.r.t. $dist^+$. 
\end{enumerate} 
We will show that the existence of such a function in Appendix. 
 Let us note that the proof of its existence  essentially repeats the proof of the existence, for arbitrary $\varepsilon_1, \varepsilon_2> 0$,   of  an $\varepsilon_1$- approximation of an
  1-Lipschitz  function  by a smooth $1+\varepsilon_2$-Lipschitz function on the standard $\mathbb{R}^n$ with the standard metric.

  Now let us take the function $\frac{1}{2} \tilde f$ and consider 
  $F+ \tfrac{1}{2}d\tilde f$. This is a Finsler metric. Indeed, we need to check that $F(x,v)  + \tfrac{1}{2}d_x\tilde f(v)>0$ for all $x $ and for all $v\ne 0$. In a local coordinate system in a neighborhood of $x$ we  consider the curve 
  $t \mapsto  x + t\cdot v$, $t\in [0, \varepsilon]$.  From the definition \eqref{1} it follows that   
\begin{equation} \label{ff} F(x,v)= \lim_{t\to 0+} \tfrac{1}{t} dist^+(x, x+ t v).\end{equation}
  Now, $$-d_x\tilde f(v)= \lim_{t\to 0+} \tfrac{1}{t} (\tilde f(x) - \tilde f(x +tv)) \stackrel{\eqref{lip}}{\le} \lim_{t\to 0+} \tfrac{1.5}{t} dist^+(x, x + tv) \stackrel{\eqref{ff}}{=} 1.5 F(x,v).$$
  Then, since $F(x,v)>0$ for all $v \ne 0$  we obtain $-\tfrac{1}{2}d_x\tilde f(v) <F(x,v)$  for all $v\ne 0$ 
  implying $F+ \tfrac{1}{2} d\tilde f$ is a Finsler metric.

  Let us now prove that the Finsler metric $F+ \tfrac{1}{2} d\tilde f$ is  forward and 
  backward 
  complete.  It is sufficient to show that for every $r\in \mathbb{R}_{\ge 0} $  the balls 
 $$
 B^{\pm}_r(p):= \{ x \in M \mid dist^\pm_{F + \tfrac{1}{2}d\tilde f}(p,x)\le r\}$$ 
 are compact. Indeed,  any forward-Cauchy sequence lies in $B^+_r(p)$ for sufficiently large $r$. Then, if such balls are compact, they are complete implying  our forward-Cauchy sequence converges. Similar arguments show that if all  balls $B^-_r(p) $ are compact then the metric is backward-complete. 
 
 Now, the  $\left(F + \tfrac{1}{2}d\tilde f\right)$-distance is given  by 
 \begin{eqnarray*}  dist^+_{F + \tfrac{1}{2}d\tilde f}(p,x) &=& dist^+_F(p, x)+ \frac{1}{2} (\tilde f(x) - \tilde f(p))\\&
 \ge& D^+(x) + \tfrac{1}{2}(f(x) - f(p))   -1 \\ &=& D^+(x) + \tfrac{1}{4} \left(D^-(x) - D^+(x)\right) - 1\\& =& 
  \tfrac{3}{4}D^+(x)+ \tfrac{1}{4} D^-(x)-1.\end{eqnarray*} 
  As we explained in the begining of the proof,   since the function $D^+ + D^-$ is proper, the function $\tfrac{3}{4}D^+(x)+ \tfrac{1}{4} D^-(x) - 1$ is proper as well  implying the function $dist_{F + \tfrac{1}{2}d\tilde f}^+(p, x)$ is proper implying the balls $B_r(p)$ are  compact so the metric $F + \tfrac{1}{2}d\tilde f$ is forward-complete. The proof that the metric is backward-complete is similar. Theorem \ref{thm1}  is proved.

\begin{remark} \label{rem:1}  In the proof we constructed a \underline{smooth}  function $\tilde f$ such that the trivial projective change $F+ \tfrac{1}{2} \tilde f$ is  a forward and backward complete Finsler metric. If we do not require the smoothness, we can simply take the trivial  projective change corresponding to the 
 function $\tfrac{1}{2} f$. As we have shown above, the function is locally 
 Lipschitz so its differential is defined at almost every point so the function $F + \tfrac{1}{2} df $ is defined almost everywhere. Moreover,  for every curve $c$ the formula \eqref{2} gives us a well-defined length  (because the restriction of a locally Lipschitz function to a smooth curve is locally Lipschitz) and the length in  $F + \tfrac{1}{2} df $  is related to the length in $F$ by the formula \eqref{3} so  the (not everywhere defined) Finsler metric $F + \tfrac{1}{2} df $ generates a  forward- and backward complete distance function.

\end{remark}

\section{Appendix: approximating  a Lipschitz function  by a smooth function}

Let $(M, F)$ be a Finsler manifold.  Assume the function $f$ is 1-Lipschitz w.r.t. to the distance function generated by $F$, that is for every $x,y \in M$ we have \begin{equation} dist^+(x,y)\ge f(x)- f(y). \label{lip1}\end{equation} 

Our goal is to show that 

{\it for every $\varepsilon_1, \varepsilon_2>0$ there exists a smooth function $\tilde f$ such that  \begin{itemize} \item 
$|\tilde f (x) - f(x)|<\varepsilon_1$ for all $x$ 
and  such that \item $\tilde f$ is $(1+ \varepsilon_2)$-Lipschitz w.r.t. to the distance function generated by $F$. \end{itemize} }

The special cases of this statement, when the $F$ generates an euclidean distance or is a Riemannian metrics, are known: in the euclidean case, this is a well known folklore, and in the  Riemannian  case it was proved for example in \cite{riem}. 

Let us first do it in a small neighborhood of  an arbitrary point $p$.  We assume that the closure  of the neighborhood is compact. We identify the neighborhood with a domain $U'\subseteq \mathbb{R}^n$, take a small number $r>0$ and a infinitely smooth  positive 
 function $\sigma:\mathbb{R}^n\to \mathbb{R} $ such that   its support lies in the Eucledian  $r-$ball, such that it  is spherically symmetric  with respect to $0$ and such that  the integral $\int_{\mathbb{R}^n}\sigma(\xi) d\xi = 1$.  We denote by $U$ the  interior of the set of the points $\{x \in U\mid \bar B_r(x)\in U\}$. The set  $U$ is open and, if $r$ is sufficiently  small, contains  
 the point $p$. 
 
 Now, denote by $\tilde f_p$ 
  the convolution of the function $f$  with   the function $\sigma$: 
  
  $$\tilde f_p(x)= \int_{\mathbb{R}^n} \sigma(x - \xi) f(\xi)d\xi .$$  
The function $\tilde f_p$ is defined for all  $x\in U$ and is smooth. Let us show that, if $r$ is small enough, 
$ |\tilde f (x) - f(x)|\le \varepsilon_1$ for all $x$ 
and $\tilde f$ is $(1+ \varepsilon_2)$-Lipschitz.

Since the function    $f$ is Lipschitz, it  
 is uniformly continuous on $U$ so for sufficiently small $r$ we have  
$|f(x)- f(y) | <\varepsilon_1$ for all $x,y \in U $ such that $d(x,y)$ .   
We consider \begin{eqnarray*} |\tilde f_p (x) - f(x)| &=& \left| \int_{\mathbb{R}^n} \sigma(x - \xi) f(\xi)d\xi  - \int_{\mathbb{R}^n} \sigma(x - \xi) f(x)d\xi\right|\\& =&  \left| \int_{\mathbb{R}^n} \sigma(x - \xi) (f(\xi) - f(x))d\xi\right|
\\ &\le&  \int_{\mathbb{R}^n} \sigma(x - \xi) \varepsilon_1 d\xi = \varepsilon_1.\end{eqnarray*}

Let us show that, for a sufficiently small $r$, the function $\tilde f_p$ is $1$-Lipschitz. Since the function is smooth, it is sufficient to show that for every  $x $ and for every $v$ 
 the directional  derivative of the function $\tilde f_p$  at the point $x$ in the direction $v$ 
  is less than $(1+ \varepsilon_2)F(x,v)$. Without loss of generality we can think that $v= \frac{\partial}{\partial {x^1}}$

We have 
$$\frac{\partial}{\partial {x^1}} \int_{\mathbb{R}^n}  \sigma(x -\xi)  f(\xi)d\xi \stackrel{\eqref{i1}}{ =}   \int_{\mathbb{R}^n} f(\xi)\frac{\partial}{\partial {x^1}} \sigma(x -\xi)  d\xi \stackrel{\eqref{i2}}{ =}  -\int_{\mathbb{R}^n} f(\xi)\frac{\partial}{\partial {\xi^1}} \sigma(x -\xi)  d\xi$$ $$ \stackrel{\eqref{i3}}{ =}   \int_{\mathbb{R}^n}  \sigma(x -\xi)  df\stackrel{\eqref{i4}}{ \le} \int_{\mathbb{R}^n}  \sigma(x -\xi)  F\left(\xi, \tfrac{\partial}{\partial {x^1}}\right)d\xi.  $$

Let us explain the equalities/inequalities  in the formula above. \begin{enumerate} \item \label{i1} Here we used the standard formula of the differentiation of an integral depending on the parameter (in this case, $x^1$). 
\item  \label{i2} Here we used that $\frac{\partial}{\partial {x^1}} \sigma(x - \xi)= - \frac{\partial}{\partial {\xi^1}}\sigma(x - \xi)$. \item  \label{i3} Here we used  integration by parts: $\int udv = uv| - \int vdu $. The role of the function $v$ plays the function $\sigma(x -\xi)$. The role of of  $u$ plays the functions $f(x)$  considered as a function of one variable $x^1$. 
Since $f$ is Lipschitz, is has bounded variation so $df$ is a well defined measure. 

Now, for this choice of the functions $u, v$,  the term $uv|$ disappears since  the function  $u= \sigma(x- \xi) $ has compact support, so we obtain   $\int udv = -\int vdu$ which gives us \eqref{i3}.
\item \label{i4}  Here we used that $\sigma$ is nonnegative and that the measure $F\left(\xi,\tfrac{\partial}{\partial {x^1}}\right) d\xi$ is greater than $ df$.  

\end{enumerate} 

Now, since the Finsler function $F$ is continuous, it is uniformally continuous  on the unite spherical
 bundle $S_1U$ implying that, if $r$ is sufficiently small, $F\left(\xi,\tfrac{\partial}{\partial {x^1}}\right)$  is $\varepsilon_2\cdot F\left(x,\tfrac{\partial}{\partial {x^1}}\right)-$close to $F\left(x,\tfrac{\partial}{\partial {x^1}}\right)$ for $\xi$ that are $r-$close to $x$. Then, $$\int_{\mathbb{R}^n}  \sigma(x -\xi)  F\left(\xi, \tfrac{\partial}{\partial {x^1}}\right)d\xi \le \int_{\mathbb{R}^n}  \sigma(x -\xi)  (1 + \varepsilon_2) F\left(x, \tfrac{\partial}{\partial {x^1}}\right)d\xi = (1 + \varepsilon_2) F\left(x, \tfrac{\partial}{\partial {x^1}}\right)$$
implying that the $v$-derivative of the 
  function $\tilde f$ is not greater than  $(1 + \varepsilon_2) F(x, v)$ implying the function $\tilde f$ is $(1 + \varepsilon_2)$-Lipschitz.

Thus, for every point $p$  we can choose a neighborhood  $U_p$ such that for every $\tilde \epsilon_1>0, \tilde \varepsilon_2>0$ 
we can $\tilde \varepsilon_1$-approximate $f$  in the neighborhood 
$U_p$ by a $(1+ \tilde \varepsilon_2)$-Lipschitz function $\tilde f_p$ on  $U_p$. 
We take a locally finite cover  $U_p$, $p\in P$  of 
$M$ by such neighborhood and choose a smooth partition of unity  $\mu_p$ corresponding to this cover. We think that  the  approximations functions $\tilde f_p$ are defined on the whole manifold  (though  it is not important what values do the functions $\tilde f_p$  have on the points that do not lie in $U_p$ since in all formulas below  we will multiply $\tilde f_p$ by $\mu_p$ and   all $\mu_p$ are zero outside of $U_p$).

  Now, set 
\begin{equation}\label{deff} \tilde f:= \sum_{p}\tilde f_p \cdot \mu_p.\end{equation}   The function $\tilde f$ is well-defined since in a small neighborhood of  every point $x$ 
only finite  many terms of the sum  are not zero, and  is evidently smooth.  
Let us  show that we can chose the numbers $\varepsilon_1(p), \varepsilon_2(p)$  for every  $U_p$ such that 
the function $\tilde f$ satisfies  our requirements.

Suppose  a point  $x$ lies in the intersection of $k$ neighborhoods of the cover $U_p$, which we denote by 
$U_1,...,U_k$. We will denote by $\mu_1,...,\mu_k$ the corresponding elements of the partition of unity and by $\tilde f_1,...,\tilde f_k$ the correspondent approximation $\tilde f_p$; we will show that there exists $\tilde \varepsilon_1, \tilde\varepsilon_2>0$ such that, if $f_i$ are 
  $(1+ \tilde \varepsilon_1)$-approximations of the restriction of the functions $f$ to $U_i$, then  $\tilde f$ is a 
   $(1+ \varepsilon_1)$-approximations of the restriction of the functions $f$ to $\bigcap_{i=1}^kU_i$. 
   
Indeed,

$$\tilde f(x)- f(x) = \mu_1(x) (\tilde f_1(x) - f(x))+ ... + \mu_k (\tilde f_k(x) - f(x))\le k\cdot 
\tilde \varepsilon_1.$$
Thus, for a $\tilde \varepsilon_1<\tfrac{1}{k} \varepsilon_1$ the function $\tilde f$ is   indeed  an  $\varepsilon_1$-approximation of  $f$. 

Let us now show that for a sufficiently small $\tilde\varepsilon_1, \tilde\varepsilon_2>0 $ the function $\tilde f$  is indeed  $(1+ \varepsilon_2)$-Lipschitz. 

It is sufficient to prove that for every tangent vector $v$ 
the directional derivative  of $\tilde f$ 
 in the direction $v$ is less than or equal to $(1+ \varepsilon_2)F(x,v)$. We take a point 
 $x$ such that it lies in the intersection of $k$ elements of the cover  which we  again denote by $U_1,...,U_k$. 
 Then $\tilde f$ given by  \eqref{deff} is actually a finite sum  $$
 \tilde f= \mu_1\tilde f_1 + \cdots + \mu_k \tilde f_k.$$  Without loss of generality we can think that $v= \tfrac{\partial }{\partial x^1}$; we need to show that 
 $$\tfrac{\partial }{\partial x^1} (\mu_1 \tilde f_1 + \cdots + \mu_k\tilde f_k) \le (1+ \varepsilon_2)F\left(x,\tfrac{\partial }{\partial x^1}\right).$$
The left hand side of the above inequality is equal, in view of the equalities $\tfrac{\partial }{\partial x^1}\mu_1 + \cdots + \tfrac{\partial }{\partial x^1} \mu_k= \tfrac{\partial}{\partial x^1}1 =  0$, to 
$$\begin{array}{ll} &  \tfrac{\partial }{\partial x^1} (\mu_1 \tilde f_1 + \cdots + \mu_k\tilde f_k) - 
  \tilde f \cdot \left(  \tfrac{\partial }{\partial x^1}\mu_1 + \cdots + \tfrac{\partial }{\partial x^1} \mu_k \right) \\ {=}& 
(\tilde f_1 - \tilde f)\tfrac{\partial }{\partial x^1}\mu_1 +  \cdots+ (\tilde f_{k} - \tilde f)\tfrac{\partial }{\partial x^1}\mu_{k}\\ +& \underbrace{\mu_1\tfrac{\partial }{\partial x^1}\tilde f_1 + \cdots+  \mu_k\tfrac{\partial }{\partial x^1}\tilde f_k}_{ \le (1+ \tilde \varepsilon_2)F\left(x,\tfrac{\partial }{\partial x^1}  \right)}.\end{array}$$
Now, since the functions $\mu_i$ have bounded support,  the derivatives  $\tfrac{\partial }{\partial x^1}\mu_i$ are bounded implying that  the sum above   
  is less than 
$(1+  \varepsilon_2)F\left(x,\tfrac{\partial }{\partial x^1}  \right)$ for sufficiently small $\tilde \varepsilon_1, \tilde \varepsilon_2$  as we claimed.

\subsection*{Acknowledgement} The work was started  during the  VI International Meeting on Lorentzian Geometry (Granada, September 6--9, 2011) and was initiated by the questions by  E. Caponio,  M.A. Javaloyes and  M. S\'anchez; I thank them for this and for the stimulating discussions at the final stage of the preparation of the paper, and  Mike Scherfner  for pointing out a misprint. 
The standard proof of the existence of a smooth Lipschitz  approximation of a Lipschitz function on $\mathbb{R}^n$ whose generalization for  Finsler  metrics is the main result of Appendix was explained to me by Yu. Burago, S. Ivanov and A. Petrunin.


\begin{thebibliography}{99}
 \bibitem{riem} 
D. Azagra, J. Ferrera, F. Lopez-Mesas, Y. Rangel,  {\it Smooth approximation of Lipschitz functions
on Riemannian manifolds,} J. Math. Anal. Appl. {\bf 326}(2007) 1370--1378.

\bibitem{R3}  A. N. Bernal, M. S\'anchez, {\it  On smooth Cauchy hypersurfaces and
Gerochís splitting theorem,}  Comm. Math. Phys. {\bf 243}(2003), no. 3, 461--470.

\bibitem{R5}  A. N. Bernal, M. S\'anchez, {\it Smoothness of time functions and the
metric splitting of globally hyperbolic spacetimes,} Comm. Math. Phys. {\bf 257}(2005), no. 1, 43--50.

 
 

\bibitem{1} E. Caponio, M. A. Javaloyes, M. S\'anchez, {\it  
On the interplay between Lorentzian Causality and Finsler metrics of Randers type,} Rev. Mat. Iberoamericana  
{\bf 27}(2011), no. 3, 919--952,  arXiv:0903.3501
 
 \bibitem{2} E. Caponio, M. A. Javaloyes, A. Masiello, {\it On the energy functional on Finsler manifolds and applications to stationary spacetimes},  	Math. Ann., {\bf 351}(2011)  365--392,  	arXiv:math/0702323v4
 
 \bibitem{2a} E. Caponio, M.A. Javaloyes, A. Masiello, {\it Morse theory of causal geodesics in a
stationary spacetime via Morse theory of geodesics of a Finsler metric,} 
Annales de l'Institut Henri Poincare (C) Non Linear Analysis {\bf 27}(2010),  857--876, 2010. 

\bibitem{2b} E. Caponio, A. V. Germinario, M. S\'anchez, {\it  Geodesics on convex regions of stationary spacetimes and Finslerian Randers spaces}, arXiv:1112.3892

 
 \bibitem{DPS} A. Dirmeier, M. Plaue, M. Scherfner, {\it Growth conditions, Riemannian completeness and Lorentzian causality}, to appear in J. Geom. Phys. (doi:10.1016/j.geomphys.2011.04.017)
 
 


\bibitem{R1} M. A.  Javaloyes, M. S\'anchez,  {\it A note on the existence of standard
splittings for conformally stationary spacetimes. } Classical Quantum Gravity
{\bf 25}(2008), no. 16, 168001.

\bibitem{4} I. Kovner, {\it Fermat principles for arbitrary space-times}, Astrophysical Journal, {\bf 351}(1990),
pp. 114--120.

\weg{
 \bibitem{short}
V. S. Matveev,
\emph{Geometric explanation of Beltrami theorem},
Int. J. Geom. Methods Mod. Phys. \textbf{3} (2006), no. 3, 623--629.

\bibitem{CMH}
V. S. Matveev,
\emph{Lichnerowicz-Obata conjecture in dimension two,}
Comm. Math. Helv. \textbf{81}(2005) no. 3,  541--570.
\bibitem{archive}
V. S. Matveev,
\emph{Proof of projective Lichnerowicz-Obata conjecture},
 J. Diff. Geom. (2007), {\bf 75}(2007),  459--502,  arXiv:math/0407337}
 
 
 
\weg{\bibitem{5} V. Perlick, {\it On Fermat's principle in general relativity. I. The general case,} Classical Quantum
Gravity, {\bf 7}(1990), pp. 1319ñ1331}

\weg{
\bibitem{schur} F. Schur, {\it Ueber den Zusammenhang der R\"aume constanter
Riemann'schen Kr\"ummumgsmaasses mit den projektiven R\"aumen, } Math. Ann. {\bf 27}(1886), 537--567.

\bibitem{dedicata}
P. J. Topalov and V. S. Matveev,
\emph{Geodesic equivalence via integrability},
Geometriae Dedicata \textbf{96} (2003), 91--115.}
\end{thebibliography}
\end{document}